\documentclass[letterpage, 11pt, notitlepage]{article}

\usepackage[margin=1.0in]{geometry}

\usepackage{epsfig}
\usepackage{amssymb}
\usepackage{amsmath}
\usepackage{amsfonts}
\usepackage{graphicx}
\usepackage{hyperref}
\usepackage{graphicx}
\usepackage{caption}
\usepackage{subcaption}
\usepackage{color}
\usepackage{amsthm,epstopdf,titling,url,array}

\theoremstyle{plain}
\newtheorem{thm}{Theorem}
\newtheorem{lem}{Lemma}
\newtheorem{prop}{Proposition}
\newtheorem{cor}{Corollary}

\theoremstyle{definition}

\newtheorem{conj}{Conjecture}
\newtheorem{exmp}{Example}

\theoremstyle{remark}
\newtheorem*{rem}{Remark}

\begin{document}

\title{\fontsize{25}{25}\selectfont A queueing system with on-demand servers: local stability of fluid limits}
\author{
Lam M. Nguyen\\
Department of Industrial \\ and Systems Engineering\\ 
Lehigh University \\ 
Bethlehem, PA 18015\\
lmn214@lehigh.edu\\
  \and
Alexander L. Stolyar\\
Department of Industrial and Enterprise \\ Systems Engineering\\ 
University of Illinois at Urbana-Champaign\\
Urbana, IL 61801\\
stolyar@illinois.edu\\
}
\maketitle

\begin{abstract}
We study a system, where a random flow of customers is served by servers (called agents) invited on-demand. Each invited agent arrives into the system after a random time; after each service completion, an agent returns to the system or leaves it with some fixed probabilities. Customers and/or agents may be impatient, that is, while waiting in queue, they leave the system at a certain rate (which may be zero). We consider the queue-length-based feedback scheme, which controls the number of pending agent invitations, depending on the customer and agent queue lengths and their changes. The basic objective is to minimize both customer and agent waiting times.

We establish the system process fluid limits in the asymptotic regime where the customer arrival rate goes to infinity. We use the machinery of switched linear systems and common quadratic Lyapunov functions to approach the stability of fluid limits at the desired equilibrium point, and derive a variety of sufficient local stability conditions. For our model, we conjecture that local stability is in fact sufficient for global stability of fluid limits; the validity of this conjecture is supported by numerical and simulation experiments. When local stability conditions do hold, simulations show good overall performance of the scheme.

\end{abstract}


\section{Introduction}\label{intro}

Consider a service system where a random flow of customers arrive exogenously. Servers, called \textit{agents}, can be invited on-demand at any time. Invited agents arrive into the system not immediately, but after a random delay. When a customer is matched with an agent, a service occurs. After completing the service, the agent can either leave the system or return to serve more customers. Customers and/or agents may be impatient, that is, they abandon the system if their wait in queue exceeds some random \textit{patience time}. The objective is to keep waiting times of both customers and agents small. Such system is schematically shown in Figure \ref{system_abandonment}. \\

The model we consider is a generalized version of that in \cite{lam2015agents,Pang@2014}. In \cite{Pang@2014}, there is no abandonment for both queues, and agents always leave the system after service completions. The model in \cite{lam2015agents} also has no abandonment, but, like in our model, an agent may return to the system after a service completion. Thus, our model is more realistic in many scenarios because customer abandonment is a key factor for call center operations (see e.g. \cite{Garnett@2002,Zeltyn2005}). \\

More specifically, the model in this paper is as follows. Customers arrive as a Poisson process and join a customer queue if no agent is available. Agents can be invited into the system exogenously, and join an agent queue after a random exponentially distributed time. There is an infinite pool of potential agents, which can be invited to serve customers. Customer service times are i.i.d. exponential. After the service completion, the customer leaves the system while the agent can return to the agent queue with some fixed probability. The matching of customers and agents is done in first-come-first-served (FCFS) order. The head-of-the-line customer and agent are matched immediately and together go to service, that is, there cannot be non-zero number of customers and agents simultaneously in the customer and agent queues. Customers and/or agents may be impatient and the patience times are independently exponentially distributed. \\

The model is primarily motivated by call/contact centers (see \cite{stolyar2010pacing}), where agents that we consider are highly skilled. It is not reasonable to set a fixed working schedule for these agents since their time is very valuable. Instead, they are invited on-demand in real time. The purpose is to design a real-time adaptive agent invitation scheme that minimizes customer and agent waiting times. However, designing an effective, simple and robust agent invitation strategy is non-trivial due to randomness in agent behavior. \\

We study a feedback-based adaptive scheme of \cite{stolyar2010pacing,Pang@2014,lam2015agents}, called \textit{queue-length-based feedback scheme}, which controls the number of pending agent invitations, depending on the customer and/or agent queue lengths and their changes. The algorithm analysis in this paper is substantially more challenging due to greater generality of our model. Just like in \cite{lam2015agents,Pang@2014}, we consider a ``stylized'' version of the invitation scheme to make the analysis more tractable. Our simulation experiments in section \ref{fluid_actual} show that the behavior of the stylized scheme is very close to that of the more practical version of the queue-length-based feedback scheme. \\

We consider the system in the asymptotic regime where the customer arrival rate goes to infinity while the distributions of the agent response times, the service times and the patience times are fixed. We show convergence of the fluid-scaled process to the fluid limit (Theorem \ref{thrm2}), which satisfies a system of differential equations. The key property of interest is the convergence of the fluid limit trajectories to the equilibrium point (at which the queues are zero). This property is referred to as {\em global stability} of the fluid limits. Establishing global stability appears to be very challenging, due to the fact that fluid limits have complicated behavior -- there are two domains where they follow different ODEs, and a ``reflecting'' boundary. In this paper, we focus on the {\em local stability} of fluid limits, defined as the stability of the dynamic system which describes fluid limit trajectories away from the boundary. The \textbf{main results} in this paper (Theorem \ref{thrm3}) give sufficient local stability conditions; the proof uses  the machinery of switched linear systems and common quadratic Lyapunov functions \cite{lin@2009,Shorten@2007}. Theorem \ref{thrm3} implies many useful sufficient local stability conditions (Corollaries \ref{cor1} - \ref{cor10}) for special cases, including those where customers never abandon or agents certainly leave the system after service completions. (Some of these corollaries -- namely, Corollaries \ref{cor3}, \ref{cor4} and \ref{cor10} -- strengthen the results in \cite{lam2015agents}  for the non-abandonment system.) These sufficient local stability conditions are robust and easy to achieve in practice.  Finally, we conjecture that, for our model, local stability is in fact sufficient for global stability, based on a large number of numerical and simulation experiments. Our simulation experiments also show good overall performance of the feedback scheme when the local stability conditions do hold. \\

The model has many applications, or potential applications. For a general discussion of modern call/contact centers and their management, see, e.g. \cite{AAM,liveops}. Another example is telemedicine \cite{WinNT}, where ``agents'' are doctors, invited on-demand to serve patients remotely. The model also arises in other applications, such as crowdsourcing-based customer service (see e.g. \cite{Arise1,Arise2}), taxi-service system, buyers and sellers in a trading market, and assembly systems. The model has relation to classical assemble-to-order models, where customers are orders and ``invited agents'' are products, which cannot be produced/assembled instantly. The model is also related to ``double-ended queues'' (see e.g. \cite{K66,LGK}) and matching systems (see e.g. \cite{Gurvich@2014}); although in such models arrivals of all types into the system are typically exogenous, as opposed to being controlled. \\

\textbf{Paper organization}. The rest of the paper is organized as follows. Some  background facts on switched linear systems and common quadratic Lyapunov functions are given in section \ref{necessaryfact}. In section \ref{model}, we describe the model and algorithm in detail. Section \ref{mainresults} states the main results of the paper, which are proved in sections \ref{proof2} and \ref{proof3}. Section \ref{numerical} provides numerical and simulation experiments; it also contains our conjectures about global and local stability of fluid limits, supported by these experiments. A discussion of the results and future work is in section \ref{conclusion}. \\

\textbf{Basic notation}: Symbols $\mathbb{N}$, $\mathbb{Z}$, $\mathbb{R}$, $\mathbb{R}_{+}$ denote the sets of natural, integer, real, real non-negative numbers, respectively. $\mathbb{R}^d$ denotes the $d$-dimensional vector space. $\mathbb{R}^{d \times d}$ denotes the set of all $d \times d$ real matrices. The standard Euclidean norm of a vector $x \in \mathbb{R}^n$ is denoted $\|x\|$. For a vector $a$ and matrix $A$, we write their transposes as $a^T$ and $A^T$, respectively. For a matrix $A$, we write its inverse and determinant as $A^{-1}$ and $\text{det}(A)$, respectively. We write $x(\cdot)$ to mean the function (or random process) $(x(t), t \geq 0)$. For a real-valued function $x(\cdot): \mathbb{R}_{+} \to \mathbb{R}$, we use either $x^{\prime}(t)$ or $(d/dt)x(t)$ to denote the derivative, and for $x(\cdot): \mathbb{R}_{+} \to \mathbb{R}^d$, $(d/dt)x(t) = (x^{\prime}_1(t),\dots,x^{\prime}_d(t))$. For $x \in \mathbb{R}$, $x^{+} = \max\{x,0\}$ and $x^{-} = - \min\{x,0\}$; and $\text{sgn}(x) = 1$ if $x > 0$, $\text{sgn}(x) = 0$ if $x = 0$, and $\text{sgn}(x) = -1$ if $x < 0$.
For $x,y \in \mathbb{R}$, we denote $x \wedge y = \min\{x,y\}$ and $x \vee y = \max\{x,y\}$. $a \Leftrightarrow b$ means ``$a$ is equivalent to $b$''; $a \Rightarrow b$ means ``$a$ implies $b$''. We write $x^r \to x \in \mathbb{R}^n$ to denote ordinary convergence in $\mathbb{R}^n$. For a finite set of scalar functions $f_n(t)$, $t \geq 0$, $n \in \mathbb{N}$, a point $t$ is called \textit{regular} if for any subset $\mathbb{N}_0 \subseteq \mathbb{N}$, the derivatives 
\begin{gather*}
\frac{d}{dt}\max_{n \in \mathbb{N}_0} f_n(t) \ \text{and} \ \frac{d}{dt}\min_{n \in \mathbb{N}_0} f_n(t)
\end{gather*}
exist. (To be precise, we require that each derivative is proper: both left and right derivatives exist and are equal.) \\


\textbf{Abbreviations}: \textit{u.o.c.} means \textit{uniform on compact sets} convergence of functions, with the argument  determined by the context (usually in $[0,\infty)$); \textit{w.p.1} means \textit{with probability 1}; \textit{i.i.d.} means \textit{independent identically distributed}; RHS means \textit{right hand side}; FSLLN means \textit{functional strong law of large numbers}; \textit{CQLF} means \textit{common quadratic Lyapunov function}; \textit{LTI system} means \textit{linear time-invariant system}.

\section{Switched linear systems and CQLF}\label{necessaryfact}

Common quadratic Lyapunov functions for switched linear systems play an important role in deriving our results. In this section, we provide some necessary background. \\

Consider a \textit{switched linear system}
\begin{gather}\label{swichedsystem}
\Sigma_S: u^\prime(t) = A(t) u(t) \ , \ A(t) \in \mathcal{A} = \{A_1, \dots, A_m\}
\end{gather}

where $\mathcal{A}$ is a set of matrices in $\mathbb{R}^{n \times n}$, and $t \to A(t)$ is a 
mapping from nonnegative real numbers into $\mathcal{A}$. (Usually, as in \cite{Shorten@2007}, this
mapping is required to be piecewise constant with only finitely many discontinuities in any bounded time-interval.
In our case this additional condition is not important, because our switched system will have a continuous derivative; see equation (\ref{noboundary_system2}) below.) For $1 \leq i \leq m$, the $i^{th}$ constituent system of the switched linear system (\ref{swichedsystem}) is the \textit{linear time-invariant (LTI) system}
\begin{gather}\label{ltisystem}
\Sigma_{A_i}: u^\prime(t) = A_i u(t). 
\end{gather}

The origin is an \textit{exponentially stable equilibrium} of the switched linear system $\Sigma_s$ if there exist real constants $C > 0$, $a > 0$ such that $\|u(t)\| \leq C e^{-a t} \|u(0)\|$ for $t \geq 0$, for all solutions $u(t)$ of the system (\ref{swichedsystem}) (see \cite{Hespanha@2004,Shorten@2007}). \\

A symmetric square $n \times n$ matrix $M$ with real coefficients is \textit{positive definite} if $z^T M z > 0$ for every non-zero column vector $z \in \mathbb{R}^n$. A symmetric square $n \times n$ matrix $M$ with real coefficients is \textit{negative definite} if $z^T M z < 0$ for every non-zero column vector $z \in \mathbb{R}^n$. A square matrix $A$ is called a \textit{Hurwitz matrix} (or \textit{stable matrix}) if every eigenvalue of $A$ has strictly negative real part. The following fact is the Hurwitz criterion of matrices in $\mathbb{R}^{3 \times 3}$ (see \cite{pontryagin@1962}). 

\begin{prop}[\cite{pontryagin@1962}] 
\label{prop2}
Let $L(\lambda) = \det(A - \lambda I) = 0$ be the characteristic equation of matrix $A$ in $\mathbb{R}^{3 \times 3}$:
\begin{gather}
L(\lambda) = a_0 \lambda^3 + a_1 \lambda^2 + a_2 \lambda + a_3 = 0 \ , \ a_0 > 0. 
\end{gather}

Matrix $A$ is Hurwitz if and only if $a_1$, $a_2$, $a_3$ are positive and $a_1 a_2 > a_0 a_3$. 
\end{prop}

The function $V(u) = u^T P u$ is a \textit{quadratic Lyapunov function} (QLF) for the system $\Sigma_{A} : u^\prime(t) = A u(t)$ if (i) $P$ is symmetric and positive definite, and (ii) $P A + A^T P$ is negative definite. Let $\{A_1,\dots,A_m\}$ be a collection of $n \times n$ Hurwitz matrices, with associated stable LTI systems $\Sigma_{A_1},\dots,\Sigma_{A_m}$. Then the function $V(u) = u^T P u$ is a \textit{common quadratic Lyapunov function} (CQLF) for these systems if $V$ is a QLF for each individual system (see \cite{lin@2009,Shorten@2007}). \\

The following facts will be used in the proof of our main results (Theorem \ref{thrm3}). 

\begin{prop}[\cite{lin@2009,Shorten@2007}]
\label{prop1}
The existence of a CQLF for the LTI systems is sufficient for the exponential stability of the switched linear system. 
\end{prop}

\begin{prop}[\cite{lin@2009,Shorten@2007}]
\label{prop4}
Let $A_1$ and $A_2$ be Hurwitz matrices in $\mathbb{R}^{n \times n}$, and the difference $A_1 - A_2$ has rank one. Then two systems $u^\prime(t) = A_1 u(t)$ and $u^\prime(t) = A_2 u(t)$  have a CQLF if and only if the matrix product $A_1 A_2$ has no negative real eigenvalues. 
\end{prop}

\begin{prop}[\cite{shorten@2004}]
\label{prop5}
If $A_1^{-1}$ is non-singular, the product $A_1 A_2$ has no negative eigenvalues if and only if $A_1^{-1} + \tau A_2$ is non-singular for all $\tau \geq 0$. 
\end{prop}

\section{Model and algorithm}\label{model}

\subsection{Model}

\begin{figure}[h]
  \begin{center}
   \includegraphics[width=0.6\textwidth]{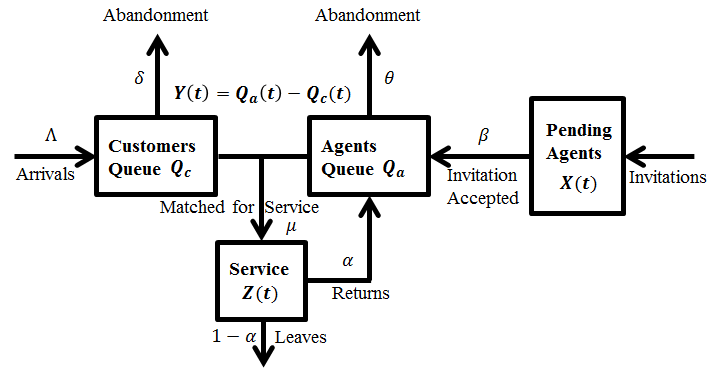}
   \caption{An agent invitation system}
   \label{system_abandonment}
  \end{center}  
\end{figure}

Our model is a generalization of that considered in \cite{lam2015agents,Pang@2014}.
Customers arrive according to a Poisson process of rate $\Lambda > 0$, and join a customer queue waiting for an available agent and are served in the order of their arrival. There is an infinite pool of 'potential' agents, which can be invited to serve customers. After a potential agent is invited, it becomes a 'pending' agent; we refer to such an event as an {\em invitation}. A pending agent 'accepts' its invitation and becomes 'active' agent after a random, exponentially distributed, time with mean $1/\beta$; we refer to such an event as an {\em acceptance}.
 Upon acceptance events, the new active agents join the (active) agent queue. The customer and agent queues cannot be positive simultaneously: the head-of-the-line customer and agent are immediately matched, leave their queues, and together go to service. Each service time is an exponentially distributed random variable with mean $1/\mu$; after the service completion, the customer leaves the system, while the corresponding agent either remains active and rejoins the agent queue -- this occurs with probability $\alpha\in[0,1)$ -- or leaves the system with probability $1-\alpha$.  Thus, there are two ways in which agents join the queue -- when an agent becomes active (upon acceptance event) and already active agents rejoining the queue after service completions. The patience times of customers and agents are independent sequences of i.i.d. exponential random variables with rate $\delta \geq 0$ and $\theta \geq 0$, respectively. When its patience time expires while a customer or server wait in queue, they leave the system. (The model in \cite{lam2015agents} is a special case of ours, with $\delta = 0$ and $\theta = 0$; in other words, customers and agents certainly wait in their queues until they are matched. The model in \cite{Pang@2014} is a special case of ours, with $\delta = 0$, $\theta = 0$ and $\alpha = 0$.) Figure \ref{system_abandonment} depicts such a system. \\

Let $X(t)$ be the number of pending agents at time $t$. Let $Y(t) = Q_a(t) - Q_c(t)$ be the difference between the agent and customer queue lengths at time $t$. (Note that $Q_a(t) = Y^{+}(t)$ and $Q_c(t) = Y^{-}(t)$.) Let $Z(t)$ be the number of customers (or agents) in service at time $t$. The system state at time $t$ is $(X(t),Y(t),Z(t))$. 

\subsection{Algorithm}

The queue-length-based feedback scheme in \cite{lam2015agents,Pang@2014,stolyar2010pacing}, referred to as the \textit{actual scheme}, maintains a ``target'' $X_{target}(t)$ for the number of pending agents $X(t)$. $X_{target}(t)$ is changed by $\Delta X_{target}(t) = [-\gamma \Delta Y(t) - \epsilon Y(t) \Delta t]$ at each time $t$ when $Y(t)$ changes by $\Delta Y(t)$ ($+1$ or $-1$), where $\gamma > 0$ and $\epsilon > 0$ are the algorithm parameters and $\Delta t$ is the time duration from the previous change of $Y$. New agent invitations occur (i.e., the number of pending agents increases) if and only if $X(t) < X_{target}(t)$, where $X(t)$ is the actual number of pending agents; therefore, $X(t) \geq X_{target}(t)$ holds at all times. In addition, $X_{target}(t) \geq 0$; i.e. if an update of $X_{target}(t)$ makes it negative, its value is immediately reset to zero. Note that $X_{target}(t)$ is not necessarily an integer. \\

Just like in \cite{lam2015agents,Pang@2014}, to simplify our theoretical analysis, we consider a ``stylized'' version of the actual scheme, referred to as the \textit{stylized scheme}, which has the same basic dynamics, but keeps $X_{target}(t)$ integer and assumes that $X(t) = X_{target}(t)$ at all times; the latter is equivalent to assuming that not only agents can be invited instantly, but pending agents can be removed from the system at any time. Formally, the stylized scheme is defined as follows. There are six types of mutually independent, and independent of the past, events that affect the dynamics of $X(t)$, $Y(t)$ and $Z(t)$ in a small time interval $[t, t + dt]$: 
\begin{itemize}
\item a customer arrival with probability $\Lambda dt + o(dt)$,
\item an acceptance with probability $\beta X(t) dt + o(dt)$, 
\item an additional event (we will call it a type-3 event) with probability $\epsilon |Y(t)| dt + o(dt)$; unlike other events, it is triggered by the algorithm itself, as opposed to other events triggered by customers' and/or agents' ``movement'' in the system,
\item a service completion with probability $\mu Z(t) dt + o(dt)$,
\item an abandonment in the customer queue with probability $\delta Y^{-}(t) dt + o(dt)$, 
\item an abandonment in the agent queue with probability $\theta Y^{+}(t) dt + o(dt)$. 
\end{itemize}


The changes at these event times are described as follows: 
\begin{itemize}
\item Upon a customer arrival, if $Y(t) > 0$, $Z(t)$ changes by $\Delta Z(t) = 1$; and if $Y(t) \leq 0$, $Z(t)$ changes by $\Delta Z(t) = 0$. $Y(t)$ changes by $\Delta Y(t) = -1$, and $X(t)$ changes by a random quantity with average $\gamma > 0$. For example, if $\gamma = 1.7$ and $\Delta Y(t) = -1$, then $\Delta X(t) = 2$ with probability $0.7$ and $\Delta X(t) = 1$ with probability $0.3$. Note that if $\gamma$ is integer, 
$\Delta X(t) = \gamma$ w.p.1. To simplify the exposition, we assume that $\gamma > 0$ is an integer.
\item Upon an acceptance event, if $Y(t) < 0$, $Z(t)$ changes by $\Delta Z(t) = 1$; and if $Y(t) \geq 0$, $Z(t)$ changes by $\Delta Z(t) = 0$. $Y(t)$ changes by $\Delta Y(t) = 1$, and $X(t)$ changes by $\Delta X(t) = -(\gamma \wedge X(t))$, that is, the change is by $-\gamma$ but $X(t)$ is kept to be nonnegative. 
\item Upon a type-3 event, if $X(t) \geq 1$, the change $\Delta X(t) = -\text{sgn}(Y(t))$ occurs; and if $X(t) = 0$, the change $\Delta X(t) = 1$ occurs if $Y(t) < 0$ and $\Delta X(t) = 0$ if $Y(t) \geq 0$. 
\item Upon a service completion, (a) if the agent returns to the agent queue (with probability $\alpha$), then if $Y(t) < 0$, the change $\Delta Z(t) = 0$ occurs; and if $Y(t) \geq 0$, the change $\Delta Z(t) = -1$ occurs; $Y(t)$ changes by $\Delta Y(t) = 1$, and $\Delta X(t) = -(\gamma \wedge X(t))$. (b) If the agent leaves the system (with probability $1 - \alpha$), then $Z(t)$ changes by $\Delta Z(t) = -1$. 
\item Upon a customer abandonment, $Y(t)$ changes by $\Delta Y(t) = 1$, and $X(t)$ changes by $\Delta X(t) = -(\gamma \wedge X(t))$. 
\item Upon an agent abandonment, $Y(t)$ changes by $\Delta Y(t) = -1$, and $X(t)$ changes by $\Delta X(t) = \gamma$.
\end{itemize}

Let $V(t) = Y^{+}(t) + Z(t)$ be the total number of agents in the system at time $t$. 
Obviously, $(X(t),Y(t),V(t))$ is a random process with states being 3-dimensional integer vectors.
However, {\em very informally}, the basic dynamics of $(X(t),Y(t),V(t))$ under the stylized scheme can be thought of as 
described by the following ODE
\begin{gather}
\label{eq-ODE}
\begin{cases}
(d/dt)X = -\gamma(d/dt)Y - \epsilon Y \\
(d/dt)Y = \beta X - \Lambda + \alpha\mu Z + \delta Y^{-} - \theta Y^{+} \\
(d/dt)V = \beta X - (1 - \alpha) \mu Z - \theta Y^{+}.   
\end{cases}
\end{gather}
ODE (\ref{eq-ODE}) is only to provide the basic intuition for the system dynamics -- it is not used in the analysis.

\section{Main results}\label{mainresults}

We consider a sequence of systems, indexed by a scaling parameter $r \to \infty$. In the system with index $r$, the arrival rate is $\Lambda = \lambda r$, while the parameters $\alpha$, $\beta$, $\mu$, $\delta$, $\theta$, $\epsilon$, $\gamma$ do not depend on $r$. The corresponding process is $(X^r(t), Y^r(t), Z^r(t)), t\ge 0$. The desired system operating point, at which $(X^r(t), Y^r(t), Z^r(t))$ should be centered is given by $(\lambda r(1 - \alpha)/\beta, 0, \lambda r/\mu)$. 
The explanation of this choice is as follows. If an invitation scheme works as desired, $Y^r(t)$ should be close to $0$; the number of customer-agent pairs $Z^r(t)$ should be close to its average value, which is $\lambda r /\mu$, so that the customers leave the system at rate $\lambda r$; finally, $X^r(t)$ should be close to the value $\chi$, such that the total average rate at which agents join the agent queue, which is $\chi \beta + [(\lambda r)/\mu] \mu \alpha$,
is equal to the customer arrival rate $\lambda r$ -- this gives $\chi = \lambda r(1 - \alpha)/\beta$.
However, instead of considering process $(X^r(t), Y^r(t), Z^r(t))$, we will consider process $(X^r(t), Y^r(t), V^r(t))$, which is more convenient for the analysis. (Recall that $Z^r(t) = V^r(t) - (Y^{r}(t))^{+}$.) Then the natural centering value for $V^r(t)$ is same as for $Z^r(t)$, namely $\lambda r/\mu$.
We define fluid-scaled process with centering as
\begin{gather}
\label{xyv_centered} (\bar{X}^r(t),\bar{Y}^r(t),\bar{V}^r(t)) = r^{-1}\left(X^r(t) - \frac{\lambda r (1 -\alpha)}{\beta}, Y^r(t), V^r(t) - \frac{\lambda r}{\mu}\right), ~~t\ge 0. 
\end{gather}

\begin{thm} \label{thrm2}
Consider a sequence of processes $(\bar{X}^r(\cdot), \bar{Y}^r(\cdot), \bar{V}^r(\cdot))$, $r \to \infty$, with deterministic initial states such that $(\bar{X}^r(0), \bar{Y}^r(0), \bar{V}^r(0)) \to (x(0), y(0), v(0))$ for some fixed $(x(0), y(0), v(0)) \in \mathbb{R}^3$, $x(0) \geq -\frac{\lambda (1 -\alpha)}{\beta}$. Then, these processes can be constructed on a common probability space, so that the following holds. W.p.1, from any subsequence of $r$, there exists a further subsequence such that
\begin{gather}
(\bar{X}^r(\cdot), \bar{Y}^r(\cdot), \bar{V}^r(\cdot)) \to (x(\cdot), y(\cdot), v(\cdot)) \ \ u.o.c. \ \ as \ \ r \to \infty
\end{gather}

where $(x(\cdot),y(\cdot),v(\cdot))$ is a locally Lipschitz trajectory such that at any regular point $t \geq 0$
\begin{gather}
\begin{cases} \label{theorem2_system} 
x^\prime(t) = \begin{cases} 
-\gamma y^\prime(t) - \epsilon y(t), \ \text{\textit{if}} \ x(t) > -\frac{\lambda (1 -\alpha)}{\beta} \\
[-\gamma y^\prime(t) - \epsilon y(t)] \vee 0, \ \text{\textit{if}} \ x(t) = -\frac{\lambda (1 -\alpha)}{\beta}
\end{cases} \\
y^\prime(t) = \beta x(t) + \alpha \mu (v(t) - y^{+}(t)) + \delta y^{-}(t) - \theta y^{+}(t) \\
v^\prime(t) = \beta x(t) - (1 - \alpha) \mu (v(t) - y^{+}(t)) - \theta y^{+}(t). 
\end{cases} 
\end{gather}
\end{thm}

A limit trajectory $(x(\cdot),y(\cdot),v(\cdot))$ specified in Theorem \ref{thrm2} will be called a \textit{fluid limit} starting from $(x(0),y(0),v(0))$. \\ 

\begin{rem}
\label{rem-fluid-eq}
Equations (\ref{theorem2_system}), which a fluid limit must satisfy, are very natural. They can be thought of as rescaled centered versions of the (informal) equations (\ref{eq-ODE}). In addition, (\ref{theorem2_system}) includes a ``reflection'' (or, ``regulation'') at the 
boundary $x=-\frac{\lambda (1 -\alpha)}{\beta}$, i.e. condition $x(t)\ge -\frac{\lambda (1 -\alpha)}{\beta}$ is ``enforced'' as all times.
This additional condition is the centered rescaled version of the condition $X^r(t) \ge 0$, which obviously must hold at all times.
\end{rem}

Consider a dynamic system $(x(t),y(t),v(t)) \in \mathbb{R}^3$: 
\begin{gather}\label{noboundary_system2}
\begin{cases}  
x^\prime(t) = -\gamma y^\prime(t) - \epsilon y(t) \\
y^\prime(t) = \beta x(t) + \alpha \mu (v(t) - y^{+}(t)) + \delta y^{-}(t) - \theta y^{+}(t) \\
v^\prime(t) = \beta x(t) - (1 - \alpha) \mu (v(t) - y^{+}(t)) - \theta y^{+}(t). 
\end{cases} 
\end{gather}

Note that the RHS of (\ref{noboundary_system2}) is continuous. This dynamic system describes the dynamics of fluid limit trajectories when the state is away from the boundary $x = -\frac{\lambda (1 -\alpha)}{\beta}$. System (\ref{noboundary_system2}) is a generalization of the system considered in \cite{lam2015agents}, referred to as a \textit{non-abandonment system}, which is a special case of ours with $\delta = 0$ and $\theta = 0$.  \\

We say that the fluid limit is \textit{globally stable} if every fluid limit trajectory converges to the equilibrium point $(0,0,0)$;
and it is \textit{locally stable} if every trajectory of the dynamic system (\ref{noboundary_system2}) converges to the equilibrium point $(0,0,0)$. Note that exponential stability of the system (\ref{noboundary_system2}) implies local stability. \\

The following theorem is the main result of this paper. It provides sufficient exponential stability conditions for the system (\ref{noboundary_system2}).  

\begin{thm}[Sufficient exponential stability conditions] \label{thrm3}
For any set of positive $\beta$, $\mu$, $\epsilon$, $\gamma$, non-negative $\delta$ and $\theta$, and $\alpha \in [0,1)$, such that either \emph{(i)}
\begin{gather}\label{stablecond}
\gamma > \max\left\{\frac{\alpha \mu - \delta}{\beta}, \sqrt{\frac{( 2 - \alpha)\epsilon\mu + \alpha\epsilon \delta}{\beta\mu}} \right\}, 
\end{gather}

or \emph{(ii)}
\begin{gather}\label{stablecond2}
\gamma > \max\left\{ \frac{\alpha \mu - \delta + \sqrt{(\alpha \mu - \delta)^2 + 4\alpha\mu^2}}{2 \beta} , \sqrt{\max\left\{\frac{\alpha\epsilon(\delta - \mu)}{\beta\mu},0\right\}} \right\} 
\end{gather}

holds, a common quadratic Lyapunov function (CQLF) of the system \emph{(\ref{noboundary_system2})} exists, and the system \emph{(\ref{noboundary_system2})} is exponentially stable. 
\end{thm}

\begin{figure}[h]
  \begin{center}
   \includegraphics[width=1.0\textwidth]{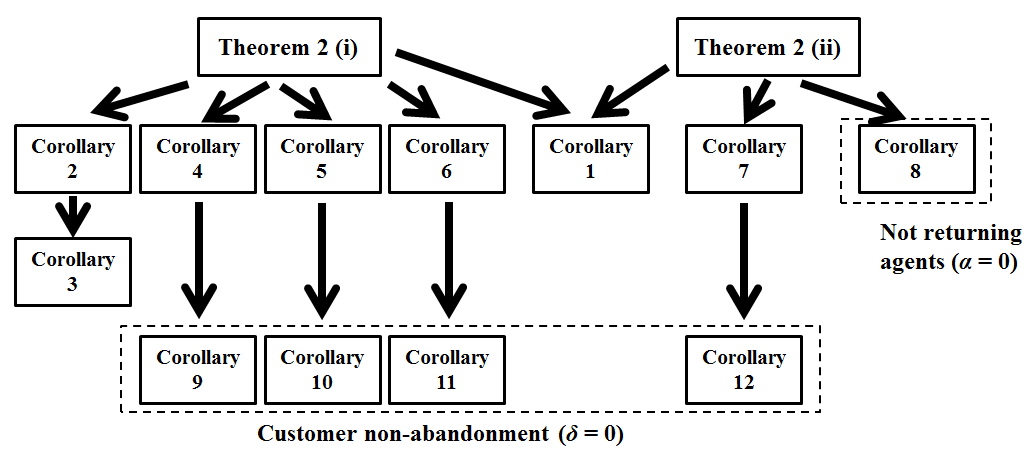}
   \caption{Result's diagram}
   \label{result_diagram}
  \end{center}  
\end{figure}

In other words, conditions (\ref{stablecond}) and (\ref{stablecond2}) are sufficient for local stability of our system. Theorem \ref{thrm3} implies the following useful sufficient local stability conditions (Corollaries \ref{cor1} - \ref{cor10}) for special cases. Figure \ref{result_diagram} depicts the connection between these results. 

\begin{cor}\label{cor1}
Given all other parameters are fixed, the system \emph{(\ref{noboundary_system2})} is exponentially stable for all sufficiently large $\gamma$. 
\end{cor}

\begin{cor}\label{cor7}
If $\alpha \mu \leq \delta$, then the system \emph{(\ref{noboundary_system2})} is exponentially stable under condition
\begin{gather}
\gamma > \sqrt{\frac{( 2 - \alpha)\epsilon\mu + \alpha\epsilon \delta}{\beta\mu}}. 
\end{gather}
\end{cor}

\begin{cor}\label{cor5}
If $\alpha \mu \leq \delta$, then the system \emph{(\ref{noboundary_system2})} is exponentially stable for all sufficiently small $\epsilon$. 
\end{cor}

\begin{cor}\label{cor8}
If $\alpha \mu > \delta$ and $\epsilon \leq \frac{(\alpha\mu - \delta)^2 \mu}{(2-\alpha)\mu\beta + \alpha\delta\beta}$, then the system \emph{(\ref{noboundary_system2})} is exponentially stable under condition
\begin{gather}
\gamma > \frac{\alpha \mu - \delta}{\beta}. 
\end{gather} 
\end{cor}

\begin{cor}\label{cor9}
If $\alpha \mu > \delta$ and $\epsilon > \frac{(\alpha\mu - \delta)^2 \mu}{(2-\alpha)\mu\beta + \alpha\delta\beta}$, then the system \emph{(\ref{noboundary_system2})} is exponentially stable under condition
\begin{gather}
\gamma > \sqrt{\frac{( 2 - \alpha)\epsilon\mu + \alpha\epsilon \delta}{\beta\mu}}. 
\end{gather} 
\end{cor}

\begin{cor}\label{cor2}
If $\alpha \mu \geq \delta$, then the system \emph{(\ref{noboundary_system2})} is exponentially stable under condition
\begin{gather}\label{oldcond}
\gamma > \frac{\alpha \mu - \delta + \sqrt{(\alpha \mu - \delta)^2 + 8 \beta \epsilon}}{2 \beta}. 
\end{gather}
\end{cor}

\begin{cor}\label{cor11}
If $\mu > \delta$, then the system \emph{(\ref{noboundary_system2})} is exponentially stable under condition
\begin{gather}
\gamma > \frac{\alpha \mu - \delta + \sqrt{(\alpha \mu - \delta)^2 + 4\alpha\mu^2}}{2 \beta}. 
\end{gather} 

(Note that this condition does not depend on $\epsilon$.)
\end{cor}

We also have the following result for the system where agents do not return to the agent queue after service completions. 
\begin{cor}\label{cor12}
If $\alpha = 0$, then the system \emph{(\ref{noboundary_system2})} is exponentially stable for all positive $\beta$, $\mu$, $\epsilon$, $\gamma$, and $\delta \geq 0$, $\theta \geq 0$. 
\end{cor}

Let us consider a special case when $\delta = 0$, referred to as a \textit{customer non-abandonment system}. Then, Corollaries \ref{cor8}, \ref{cor9}, \ref{cor2}, and \ref{cor11} imply the following sufficient local stability conditions of the customer non-abandonment system. 

\begin{cor}\label{cor3}
If $\delta = 0$, $\alpha \in (0,1)$, and $\epsilon \leq \frac{\alpha^2\mu^2}{(2-\alpha)\beta}$, then the system \emph{(\ref{noboundary_system2})} is exponentially stable under condition
\begin{gather}
\gamma > \frac{\alpha \mu}{\beta}. 
\end{gather} 
\end{cor}

\begin{cor}\label{cor4}
If $\delta = 0$, $\alpha \in (0,1)$, and $\epsilon > \frac{\alpha^2\mu^2}{(2-\alpha)\beta}$, then the system \emph{(\ref{noboundary_system2})} is exponentially stable under condition
\begin{gather}
\gamma > \sqrt{\frac{(2 - \alpha)\epsilon}{\beta}}. 
\end{gather} 
\end{cor}

\begin{cor}\label{cor6}
If $\delta = 0$, and $\alpha \in [0,1)$, then the system \emph{(\ref{noboundary_system2})} is exponentially stable under condition
\begin{gather}
\gamma > \frac{\alpha \mu + \sqrt{\alpha^2 \mu^2 + 8 \beta \epsilon}}{2 \beta}. 
\end{gather}
\end{cor}

Note that if $\theta = 0$, then Corollary \ref{cor6} is a simpler, equivalent version of the sufficient local stability condition in \cite{lam2015agents} for the non-abandonment system (Theorem 3 in \cite{lam2015agents}). Moreover, condition (\ref{stablecond2}) in Theorem \ref{thrm3} implies the following result, which does not depend on $\epsilon$, for the non-abandonment system. 

\begin{cor}\label{cor10}
If $\delta = 0$, and $\alpha \in [0,1)$, then the system \emph{(\ref{noboundary_system2})} is exponentially stable under condition
\begin{gather}
\gamma > \frac{(\alpha + \sqrt{\alpha^2 + 4\alpha})\mu}{2 \beta}. 
\end{gather}
\end{cor}

Having a variety of these sufficient local stability conditions is useful, because some or others may be easier to verify/ensure, depending on the scenario. Note that $\gamma$ and $\epsilon$ are control parameters, while all other parameters are those of the system -- they can be potentially measured/estimated in real time. It is not easy to give an intuitive meaning/interpretation of the above local stability conditions. 
Perhaps Corollary ~\ref{cor1} is the easiest to interpret: if magnitude $\gamma$ of the system response to changes in the queue length is large enough, this is sufficient for local stability.\\

\begin{rem}
We note that our {\em local} stability results apply to more general systems, exhibiting same {\em local} behavior. For example, suppose the total number of potential agents is not infinite, by finite, scaling with $r$ as $\kappa r$, where $\kappa > \lambda (1-\alpha)/\beta$. Then, the fluid limits of such system satisfy the same ODE (\ref{noboundary_system2}) in the vicinity of the origin, and therefore our local stability results apply as is. 
\end{rem}

\section{Proof of Theorem \ref{thrm2}}\label{proof2}

The proof of Theorem~\ref{thrm2} is a generalization of the proof of Theorem 1 in \cite{Pang@2014}. However, it requires additional technical details -- we present it here for completeness.\\

In order to prove Theorem \ref{thrm2}, it suffices to show that w.p.1 from any subsequence of $r$, we can choose a further subsequence, along which a u.o.c. convergence to a fluid limit holds. \\

Let $N_i(\cdot)$, $i = 1, \dots, 8$ be mutually independent unit-rate Poisson processes. $N_1$ is the process which drives customer arrivals. $N_2$ is the process which drives the acceptance of invitations. $N_3$ is the process which drives the service completions with agents leaving the system. $N_4$ is the process which drives the service completions with agents returning the agent queue. $N_5$ and $N_6$ are the processes which drive type-3 events, when variable $Y^r(t)$ is negative and positive, respectively. $N_7$ is the process which drives the abandonment of customers. $N_8$ is the process which drives the abandonment of agents. Given the initial state $(X^r(0),Y^r(0),V^r(0))$, we construct the process $(X^r(\cdot), Y^r(\cdot), V^r(\cdot))$, for all $r$, on the same probability space via a common set of independent Poisson process \cite{Pang@2007} as follows:
\begin{align} 
\label{pro_x} X^r(t) &= G^r(t) + \left( -\min_{0 \leq s \leq t} G^r(s)\right) \vee 0, \\
G^r(t) &= X^r(0) + \gamma N_1(\lambda rt) - \gamma N_2 \left(\beta \int_0^t X^r(s) ds\right) - \gamma N_4 \left(\alpha \mu \int_0^t (V^r(s) - (Y^r(s))^{+}) ds\right) - \nonumber \\
& - \gamma N_7 \left(\delta \int_0^t (Y^r(s))^{-}ds\right) + \gamma N_8 \left(\theta \int_0^t (Y^r(s))^{+}ds\right) + \nonumber \\ 
\label{pro_g} & + N_5 \left(\epsilon \int_0^t (Y^r(s))^{-}ds\right) - N_6 \left(\epsilon \int_0^t (Y^r(s))^{+}ds\right), \\
Y^r(t) & = Y^r(0) + N_2 \left(\beta \int_0^t X^r(s)ds\right) - N_1(\lambda rt) + N_4 \left(\alpha \mu \int_0^t (V^r(s) - (Y^r(s))^{+}) ds\right) + \nonumber \\ 
\label{pro_y} & + N_7 \left(\delta \int_0^t (Y^r(s))^{-}ds\right) - N_8 \left(\theta \int_0^t (Y^r(s))^{+}ds\right) , \\
V^r(t) & = V^r(0) + N_2 \left(\int_0^t \beta X^r(s) ds\right) - N_3 \left(\int_0^t (1 - \alpha) \mu (V^r(s) - (Y^r(s))^{+}) ds\right) - \nonumber \\
\label{pro_v} & - N_8 \left(\theta \int_0^t (Y^r(s))^{+}ds\right). 
\end{align}

W.p.1, for any $r$, relations (\ref{pro_x})-(\ref{pro_v}) uniquely define the realization of $(X^r(\cdot),Y^r(\cdot),V^r(\cdot))$ via the realizations of the driving processes $N_i(\cdot)$. Relation (\ref{pro_x}), the ``reflection'' at zero,  corresponds to the property that $X^r(t)$ cannot become negative. \\

The functional strong law of large numbers (FSLLN) holds for each Poisson process $N_i$: 
\begin{gather} \label{fslln}
\frac{N_i(rt)}{r} \to t \ , \ r \to \infty \ , \ \text{u.o.c.}, \ \text{w.p.1}.
\end{gather}

We consider the sequence of associated fluid-scaled processes with centering $(\bar{X}^r(\cdot),\bar{Y}^r(\cdot),\bar{V}^r(\cdot))$ as defined in (\ref{xyv_centered}). Let a constant $m > \|(x(0),y(0),v(0)\|$ be fixed. For each $r$, on the same probability space as $(\bar{X}^r(\cdot),\bar{Y}^r(\cdot),\bar{V}^r(\cdot))$, let us define a modified fluid-scaled process $(\bar{X}^r_m(\cdot),\bar{Y}^r_m(\cdot),\bar{V}^r_m(\cdot))$. Let $(\bar{X}^r_m(\cdot),\bar{Y}^r_m(\cdot),\bar{V}^r_m(\cdot))$ start from the same initial state as $(\bar{X}^r(\cdot),\bar{Y}^r(\cdot),\bar{V}^r(\cdot))$  , i.e., $(\bar{X}^r_m(0),\bar{Y}^r_m(0),\bar{V}^r_m(0)) = (\bar{X}^r(0),\bar{Y}^r(0),\bar{V}^r(0))$. The modified process  $(\bar{X}^r_m(\cdot),\bar{Y}^r_m(\cdot),\bar{V}^r_m(\cdot))$ follows the same path as $(\bar{X}^r(\cdot),\bar{Y}^r(\cdot),\bar{V}^r(\cdot))$ until the first time $t$, such that $\|(\bar{X}^r(t),\bar{Y}^r(t),\bar{V}^r(t))\| \geq m$. Denote this time by $\tau^r_m$. We then freeze the process $(\bar{X}^r_m(\cdot),\bar{Y}^r_m(\cdot),\bar{V}^r_m(\cdot))$ at the value $(\bar{X}^r(\tau^r_m),\bar{Y}^r(\tau^r_m),\bar{V}^r(\tau^r_m))$, i.e. $(\bar{X}^r_m(t),\bar{Y}^r_m(t),\bar{V}^r_m(t)) = (\bar{X}^r(\tau^r_m),\bar{Y}^r(\tau^r_m),\bar{V}^r(\tau^r_m))$ for all $t \geq \tau^r_m$.

\begin{lem}\label{lemma_xyu}
Fix $(x(0),y(0),v(0))$ and a finite constant $m > \|(x(0),y(0),v(0))\|$. Then, w.p.1  for any subsequence of $r$, there exists a further subsequence, along which $(\bar{X}^r_m, \bar{Y}^r_m, \bar{V}^r_m)$ converges u.o.c. to a Lipschitz continuous trajectory $(x_m,y_m,v_m)$, which satisfies properties \emph{(\ref{theorem2_system})} at any regular time $t \geq 0$ such that $\|(x_m(t),y_m(t),v_m(t))\| < m$. 
\end{lem} 

\textit{Proof}. For the modified fluid-scaled processes $(\bar{X}^r_m(\cdot), \bar{Y}^r_m(\cdot), \bar{V}^r_m(\cdot))$, we define the associated counting processes for upward and downward jumps. For $t \leq \tau^r_m$, 
\begin{align} 
\bar{X}^{r \uparrow}_m(t) &= r^{-1} \gamma N_1(\lambda rt) + r^{-1} \gamma N_8 \left(\theta r \int_0^t (\bar{Y}^r_m(s))^{+}ds\right) + r^{-1} N_5 \left(\epsilon r \int_0^t (\bar{Y}^r_m(s))^{-}ds\right), \\
\bar{X}^{r \downarrow}_m(t) &= r^{-1} \gamma N_2 \left(\beta r \int_0^t \left[\bar{X}^r_m(s) + \frac{\lambda (1 -\alpha)}{\beta}\right] ds\right) + r^{-1} \gamma N_4 \left(\alpha \mu r \int_0^t \left[\bar{V}^r_m(s) + \frac{\lambda}{\mu} - (\bar{Y}^r_m(s))^{+} \right] ds\right) + \nonumber \\ 
& + r^{-1} \gamma N_7 \left(\delta r \int_0^t (\bar{Y}^r_m(s))^{-}ds\right) + r^{-1} N_6 \left(\epsilon r \int_0^t (\bar{Y}^r_m(s))^{+}ds\right) , \\
\bar{Y}^{r \uparrow}_m(t) &= r^{-1} N_2 \left(\beta r \int_0^t \left[\bar{X}^r_m(s) + \frac{\lambda (1 -\alpha)}{\beta}\right] ds\right) + r^{-1} N_4 \left(\alpha \mu r \int_0^t \left[\bar{V}^r_m(s) + \frac{\lambda}{\mu} - (\bar{Y}^r_m(s))^{+} \right] ds\right) + \nonumber \\ 
& + r^{-1} N_7 \left(\delta r \int_0^t (\bar{Y}^r_m(s))^{-}ds\right)  , \\
\bar{Y}^{r \downarrow}_m(t) &= r^{-1} N_1(\lambda rt) + r^{-1} N_8 \left(\theta r \int_0^t (\bar{Y}^r_m(s))^{+}ds\right) \\
\bar{V}^{r \uparrow}_m(t) &= r^{-1} N_2 \left(\beta r \int_0^t \left[\bar{X}^r_m(s) + \frac{\lambda (1 -\alpha)}{\beta}\right] ds\right) ,  \\ 
\bar{V}^{r \downarrow}_m(t) &= r^{-1} N_3 \left((1 - \alpha) \mu r \int_0^t \left[\bar{V}^r_m(s) + \frac{\lambda}{\mu} - (\bar{Y}^r_m(s))^{+} \right] ds\right) + r^{-1} N_8 \left(\theta r \int_0^t (\bar{Y}^r_m(s))^{+}ds\right),
\end{align}

and for $t > \tau^r_m$, all these counting processes are frozen at their values at time $\tau^r_m$, that is, 
\begin{gather}
\begin{cases}
\bar{X}^{r \uparrow}_m(t) = \bar{X}^{r \uparrow}_m(\tau^r_m) \ , \ \bar{X}^{r \downarrow}_m(t) = \bar{X}^{r \downarrow}_m(\tau^r_m) \ , \\ 
\bar{Y}^{r \uparrow}_m(t) = \bar{Y}^{r \uparrow}_m(\tau^r_m) \ , \ \bar{Y}^{r \downarrow}_m(t) = \bar{Y}^{r \downarrow}_m(\tau^r_m) \ , \\
\bar{V}^{r \uparrow}_m(t) = \bar{V}^{r \uparrow}_m(\tau^r_m) \ , \ \bar{V}^{r \downarrow}_m(t) = \bar{V}^{r \downarrow}_m(\tau^r_m). 
\end{cases}
\end{gather}

Using the relations (\ref{pro_x})-(\ref{pro_v}) and the fact that for $0 \leq t \leq \tau^r_m$ the original process $(\bar{X}^r,\bar{Y}^r,\bar{V}^r)$ and the modified process $(\bar{X}^r_m,\bar{Y}^r_m,\bar{V}^r_m)$ coincide, we have for all $t \geq 0$, 
\begin{gather} 
\bar{X}^r_m(t) = \bar{G}^r_m(t) + \left(-\lambda (1 -\alpha)/ \beta - \min_{0 \leq s \leq t} \bar{G}^r_m(s)\right) \vee 0, \\
\bar{G}^r_m(t) = \bar{X}^r(0) + \bar{X}^{r \uparrow}_m(t) - \bar{X}^{r \downarrow}_m(t), \\
\bar{Y}^r_m(t) = \bar{Y}^r(0) + \bar{Y}^{r \uparrow}_m(t) - \bar{Y}^{r \downarrow}_m(t), \\
\bar{V}^r_m(t) = \bar{V}^r(0) + \bar{V}^{r \uparrow}_m(t) - \bar{V}^{r \downarrow}_m(t). 
\end{gather}

The counting processes $\bar{X}^{r \uparrow}_m(\cdot)$, $\bar{X}^{r \downarrow}_m(\cdot)$, $\bar{Y}^{r \uparrow}_m(\cdot)$, $\bar{Y}^{r \downarrow}_m(\cdot)$, $\bar{V}^{r \uparrow}_m(\cdot)$, $\bar{V}^{r \downarrow}_m(\cdot)$ are non-decreasing. Using FSLLN (\ref{fslln}) and the fact that the processes $\bar{X}^r_m(\cdot)$, $\bar{Y}^r_m(\cdot)$, and $\bar{V}^r_m(\cdot)$ are uniformly bounded by construction, we see that w.p.1. for any subsequence of $r$, there exists a further subsequence along which the set of trajectories $(\bar{X}^{r \uparrow}_m(\cdot), \bar{X}^{r \downarrow}_m(\cdot), \bar{Y}^{r \uparrow}_m(\cdot), \bar{Y}^{r \downarrow}_m(\cdot), \bar{V}^{r \uparrow}_m(\cdot), \bar{V}^{r \downarrow}_m(\cdot))$ converges u.o.c. to a set of non-decreasing  Lipschitz continuous functions $(x^{\uparrow}_m(\cdot), x^{\downarrow}_m(\cdot), y^{\uparrow}_m(\cdot), y^{\downarrow}_m(\cdot), v^{\uparrow}_m(\cdot), v^{\downarrow}_m(\cdot))$. But then the u.o.c. convergence of $(\bar{X}^r_m(\cdot), \bar{Y}^r_m(\cdot), \bar{V}^r_m(\cdot), \bar{G}^r_m(\cdot))$ to a set of Lipschitz continuous functions $(x_m(\cdot), y_m(\cdot), v_m(\cdot), g_m(\cdot))$ holds, where
\begin{gather} 
x_m(t) = g_m(t) + \left(-\lambda (1 -\alpha)/ \beta - \min_{0 \leq s \leq t} g_m(s)\right) \vee 0, \\
g_m(t) = x(0) + x^{\uparrow}_m(t) - x^{\downarrow}_m(t), \\
y_m(t) = y(0) + y^{\uparrow}_m(t) - y^{\downarrow}_m(t), \\
v_m(t) = v(0) + v^{\uparrow}_m(t) - v^{\downarrow}_m(t),
\end{gather}

and the following holds for $t$ before fluid trajectory hits $\|(x_m(t),y_m(t),v_m(t))\| = m$
\begin{align}
x^{\uparrow}_m(t) & = \gamma \lambda t + \gamma \theta \int_0^t y^{+}_m(s) ds + \epsilon \int_0^t y^{-}_m(s) ds, \\
x^{\downarrow}_m(t) & = \gamma \beta \int_0^t \left(x_m(s) + \frac{\lambda (1 -\alpha)}{\beta}\right) ds + \gamma \alpha \mu \int_0^t \left(v_m(s) + \frac{\lambda}{\mu} - y^{+}_m(s) \right) ds + \nonumber \\ & + \gamma \delta \int_0^t y^{-}_m(s) ds + \epsilon \int_0^t y^{+}_m(s) ds,  \\
y^{\uparrow}_m(t) & = \beta \int_0^t \left(x_m(s) + \frac{\lambda (1 -\alpha)}{\beta}\right) ds + \alpha \mu \int_0^t \left(v_m(s) + \frac{\lambda}{\mu} - y^{+}_m(s) \right) ds + \delta \int_0^t y^{-}_m(s) ds, \\
y^{\downarrow}_m(t) & = \lambda t + \theta \int_0^t y^{+}_m(s) ds, \\
v^{\uparrow}_m(t) & = \beta \int_0^t \left(x_m(s) + \frac{\lambda (1 -\alpha)}{\beta}\right) ds,  \\ 
v^{\downarrow}_m(t) & = (1 - \alpha)\mu \int_0^t \left(v_m(s) + \frac{\lambda}{\mu} - y^{+}_m(s) \right) ds + \theta \int_0^t y^{+}_m(s) ds.
\end{align} 

Hence,
\begin{gather}
\begin{cases}  
x^{\prime}_m(t) = \begin{cases} 
-\gamma \beta x_m(t) - \gamma \alpha \mu ( v_m(t) - y^{+}_m(t)) + \gamma \theta y^{+}_m(t) - \gamma \delta y^{-}_m(t) - \epsilon y_m(t) , \ \text{if} \ x_m(t) > -\frac{\lambda (1 -\alpha)}{\beta} \\
[-\gamma \beta x_m(t) - \gamma \alpha \mu ( v_m(t) - y^{+}_m(t)) + \gamma \theta y^{+}_m(t) - \gamma \delta y^{-}_m(t) - \epsilon y_m(t)] \vee 0, \ \text{if} \ x_m(t) = -\frac{\lambda (1 -\alpha)}{\beta}
\end{cases} \\
y{^\prime}_m(t) = \beta x_m(t) + \alpha \mu (v_m(t) - y^{+}_m(t)) + \delta y^{-}_m(t) - \theta y^{+}_m(t) \\
v{^\prime}_m(t) = \beta x_m(t) - (1 - \alpha) \mu (v_m(t) - y^{+}_m(t)) - \theta y^{+}_m(t).
\end{cases}
\end{gather}

It is easy to verify that, at any regular time $t \geq 0$ such that $\|(x_m(t),y_m(t),v_m(t))\| < m$, properties (\ref{theorem2_system}) hold for the trajectory $(x_m(\cdot),y_m(\cdot),v_m(\cdot))$. $\Box$ \\

\textit{Conclusion of the proof of Theorem \ref{thrm2}}. It is easy to see that 
\begin{gather}
\frac{d}{dt}\|(x_m(t),y_m(t),v_m(t))\| \leq C \|(x_m(t),y_m(t),v_m(t))\| \ \text{for any $m$ and some $C > 0$}. 
\end{gather}

From Gronwall's inequality \cite{opac-b1080363}, we have 
\begin{gather}
\|(x_m(t),y_m(t),v_m(t))\| \leq \|(x(0),y(0),v(0))\| e^{C t} \ \text{for all} \ t \geq 0
\end{gather}

For a given $(x(0),y(0),v(0))$, let us fix $T_l > 0$ and choose $m_l > \|(x(0),y(0),v(0)\| e^{C T_l}$. For this $T_l > 0$, there exists a subsequence $r^{l}$, along which $(\bar{X}^r, \bar{Y}^r, \bar{V}^r)$ converges uniformly to $(x_{m_l},y_{m_l},v_{m_l})$, which satisfies properties (\ref{theorem2_system}), at any $t \in [0,T_l]$. The limit trajectory $(x_{m_l},y_{m_l},v_{m_l})$ does not hit $m_l$ in $[0,T_l]$. Subsequence $r^{l} = \{r^{l}_1, r^{l}_2, \dots\}$ is such that, w.p.1, for all sufficiently large $r$ along the subsequence $r^{l}$, $(\bar{X}^r(t),\bar{Y}^r(t),\bar{V}^r(t)) = (\bar{X}^r_{m_l}(t),\bar{Y}^r_{m_l}(t),\bar{V}^r_{m_l}(t))$ at any $t \in [0,T_l]$. \\

We consider a sequence $T_1$, $T_2$, $\dots$, $\to \infty$. We construct a subsequence $r^{*}$ by using Cantor's diagonal procedure \cite{opac-b1098274} from subsequences $r^{1}$, $r^{2}$, $\dots$ ($r^{1} \supseteq r^{2} \supseteq \dots$) corresponding to $T_1$, $T_2$, $\dots$, respectively (i.e. $r^{*}_1 = r^{1}_1$, $r^{*}_2 = r^{2}_2$, $\dots$). Clearly, for this subsequence $r^{*}$, w.p.1, $(\bar{X}^r, \bar{Y}^r, \bar{V}^r)$ converges u.o.c. to $(x,y,v)$, which satisfies properties (\ref{theorem2_system}), at any regular point $t \in [0,\infty)$. $\Box$ 

\section{Proof of Theorem \ref{thrm3}}\label{proof3} 

In order to prove Theorem \ref{thrm3}, it suffices to show that LTI systems of the switched linear system (\ref{noboundary_system2}) have a CQLF. \\

The system (\ref{noboundary_system2}) is a switched linear system with $m = 2$. (Note that $y^{+} = y$ if $y \geq 0$ and $y^{+} = 0$ if $y < 0$, and $y^{-} = 0$ if $y \geq 0$ and $y^{-} = -y$ if $y < 0$.) Namely, for $y \geq 0$,  
\begin{gather}\label{system_a1}
\begin{cases}
x^\prime(t) = (-\gamma \beta) x(t) + (\gamma \alpha \mu + \gamma \theta - \epsilon) y(t) + (- \gamma \alpha \mu) v(t) \\
y^\prime(t) = (\beta) x(t) + (-\alpha \mu - \theta) y(t) + (\alpha \mu) v(t) \\
v^\prime(t) = (\beta) x(t) + ((1 - \alpha) \mu - \theta) y(t) + (-(1 - \alpha) \mu) v(t) 
\end{cases} 
\end{gather}

and for $y < 0$,  
\begin{gather}\label{system_a2}
\begin{cases}  
x^\prime(t) = (-\gamma \beta) x(t) + (\gamma \delta - \epsilon) y(t) + (-\gamma \alpha \mu) v(t) \\
y^\prime(t) = (\beta) x(t) + (-\delta) y(t) + (\alpha \mu) v(t) \\
v^\prime(t) = (\beta) x(t) + (-(1 - \alpha) \mu) v(t)
\end{cases}
\end{gather}

We can rewrite the systems above as two LTI systems $u^\prime(t) = A_1 u(t)$ and $\ u^\prime(t) = A_2 u(t)$, where $u(t) = (x(t),y(t),v(t))^T$ and
\begin{gather}\label{a1}
A_1 = \begin{pmatrix}
-\gamma \beta & \gamma \alpha \mu + \gamma \theta - \epsilon & - \gamma \alpha \mu \\
\beta & -\alpha \mu - \theta & \alpha \mu \\
\beta &  (1 - \alpha) \mu - \theta & -(1 - \alpha) \mu
\end{pmatrix} \ , \ A_2 = \begin{pmatrix}
-\gamma \beta & \gamma \delta - \epsilon & - \gamma \alpha \mu \\
\beta & -\delta & \alpha \mu \\
\beta &  0 & -(1 - \alpha) \mu
\end{pmatrix}. 
\end{gather} 

\begin{lem} \label{lemmaa1}
Matrix $A_1$ in \emph{(\ref{a1})} is Hurwitz for all positive $\beta$, $\gamma$, $\mu$, $\epsilon$, $\delta \geq 0$, $\theta \geq 0$ and $\alpha \in [0,1)$. 
\end{lem}

\textit{Proof}. The characteristic equation of $A_1$ is 
\begin{gather}
\lambda^3 + (\beta \gamma + \mu + \theta) \lambda^2 + (\beta \epsilon + \beta \gamma \mu + \mu \theta) \lambda + \beta \epsilon \mu = 0. 
\end{gather}

By Proposition \ref{prop2}, it suffices to verify that
\begin{gather}
\label{a11} \beta \gamma + \mu + \theta > 0 \ , \ \beta \epsilon + \beta \gamma \mu + \mu \theta > 0 \ , \ \beta \epsilon \mu > 0, \ \text{and} \\
\label{a12} (\beta \gamma + \mu + \theta)(\beta \epsilon + \beta \gamma \mu + \mu \theta) - \beta \epsilon \mu = \beta^2 \gamma^2 \mu + \beta^2 \gamma \epsilon + \beta \gamma \mu^2 + 2\beta\gamma\mu\theta + \beta\theta\epsilon + \mu^2\theta + \mu\theta^2 > 0. 
\end{gather} 

The conditions (\ref{a11}) and (\ref{a12}) are obviously true. $\Box$ 

\begin{lem} \label{lemmaa2}
For positive $\beta$, $\gamma$, $\mu$, $\epsilon$, $\delta \geq 0$, $\theta \geq 0$ and $\alpha \in [0,1)$, matrix $A_2$ in \emph{(\ref{a1})} is Hurwitz if and only if
\begin{gather}\label{a2hurwitz}
\left(\frac{\beta \gamma + \delta}{\mu} + (1 - \alpha)\right)\left(\frac{\beta \gamma \mu + \delta \mu (1 - \alpha)}{\beta \epsilon} + 1\right) > 1
\end{gather}
\end{lem}

\textit{Proof}. The characteristic equation of $A_2$ is 
\begin{gather}
\lambda^3 + (\beta \gamma + \mu(1 - \alpha) + \delta) \lambda^2 + (\beta \epsilon + \beta \gamma \mu + \delta \mu (1 - \alpha)) \lambda + \beta \epsilon \mu = 0. 
\end{gather}

By Proposition \ref{prop2}, it suffices to verify that
\begin{gather}\label{a21}
\beta \gamma + \mu(1 - \alpha) + \delta > 0 \ , \ \beta \epsilon + \beta \gamma \mu + \delta \mu (1 - \alpha) > 0 \ , \ \beta \epsilon \mu > 0,
\end{gather} 

and $(\beta \gamma + \mu(1 - \alpha) + \delta)(\beta \epsilon + \beta \gamma \mu + \delta \mu (1 - \alpha)) - \beta \epsilon \mu > 0$, which is equivalent to (\ref{a2hurwitz}) since
\begin{gather*}
(\beta \gamma + \mu(1 - \alpha) + \delta)(\beta \epsilon + \beta \gamma \mu + \delta \mu (1 - \alpha)) - \beta \epsilon \mu > 0 \\
\Leftrightarrow (\beta \gamma + \delta + \mu(1 - \alpha))(\beta \gamma \mu + \delta \mu (1 - \alpha) + \beta \epsilon ) > \beta \epsilon \mu \\
\Leftrightarrow \left(\frac{\beta \gamma + \delta}{\mu} + (1 - \alpha)\right)\left(\frac{\beta \gamma \mu + \delta \mu (1 - \alpha)}{\beta \epsilon} + 1\right) > 1. 
\end{gather*} 

The conditions (\ref{a21}) are obviously true. $\Box$ \\

It is easy to see that Lemma \ref{lemmaa2} implies the following result. 

\begin{cor}\label{corAhurwitz}
Matrix $A_2$ in \emph{(\ref{a1})} is Hurwitz if
\begin{gather}\label{a2hurwitz2}
\gamma > \frac{\alpha\mu - \delta}{\beta}. 
\end{gather}

(Note that $\gamma > 0$ by definition.)
\end{cor}

\begin{lem}\label{a2cond}
Matrix $A_2$ in \emph{(\ref{a1})} is Hurwitz under the condition either \emph{(\ref{stablecond})} or \emph{(\ref{stablecond2})}. 
\end{lem}

\textit{Proof}. This easily follows by applying Corollary \ref{corAhurwitz}. 

\begin{lem}\label{lemmaa1a2noeig}
Matrix product $A_1 A_2$ has no negative eigenvalues under the condition either \emph{(\ref{stablecond})} or \emph{(\ref{stablecond2})}.
\end{lem}

\textit{Proof}. With the help of MATLAB symbolic calculation, it can be shown that $A_1$ is non-singular and
\begin{gather}
A_1^{-1} = \begin{pmatrix}
-\frac{\theta}{\beta \epsilon} & -\frac{(\alpha \epsilon - \epsilon + \gamma \theta)}{\beta \epsilon} & \frac{\alpha}{\beta} \\
-\frac{1}{\epsilon} & -\frac{\gamma}{\epsilon} & 0 \\
-\frac{1}{\epsilon} & \frac{(\epsilon - \gamma \mu)}{\epsilon \mu} & -\frac{1}{\mu} \end{pmatrix}.
\end{gather}

By Proposition \ref{prop5}, to demonstrate that the product $A_1 A_2$ has no negative eigenvalues, it will suffice to show that $[A_1^{-1} + \tau A_2]$ is non-singular for all $\tau \geq 0$. We have
\begin{gather}
\det[A_1^{-1} + \tau A_2] = [\beta^2\epsilon^2\mu^2\tau^3 + \nonumber \\ + (\beta^2\epsilon^2 + \beta^2\gamma^2\mu^2 - 2\beta\epsilon\mu^2 + \delta\mu^2 \theta + \alpha\beta\epsilon\mu^2 + \beta\delta\gamma\mu^2 - \alpha\delta\mu^2\theta + \beta\gamma\mu^2\theta - \alpha\beta\delta\epsilon\mu - \alpha\beta\delta\gamma\mu^2)\tau^2 + \nonumber \\ \label{fractiondetAAA1} + (\mu^2 - \alpha\mu^2 + \beta^2\gamma^2 - 2\beta\epsilon + \delta\theta + \beta\delta\gamma + \alpha\delta\mu + \beta\gamma\theta - \alpha\mu\theta - \alpha\beta\gamma\mu)\tau + 1]/(-\beta\epsilon\mu). 
\end{gather}

To show $\det[A_1^{-1} + \tau A_2] \neq 0$ for all $\tau \geq 0$, it will suffice to show that the numerator of the ratio (\ref{fractiondetAAA1}) is strictly positive. We can represent the numerator of the ratio (\ref{fractiondetAAA1}) as follows. \\

(a) Under the condition (\ref{stablecond}), the numerator of the ratio (\ref{fractiondetAAA1}) is
\begin{gather}
\beta^2\epsilon^2\mu^2\tau^3 + (\beta\epsilon\tau - 1)^2 + \nonumber \\ + [(\beta^2\gamma^2\mu^2 - 2\beta\epsilon\mu^2 - \alpha\beta\delta\epsilon\mu + \alpha\beta\epsilon\mu^2) + \delta\mu^2 \theta(1 - \alpha) + \beta\delta\gamma\mu^2(1 - \alpha)  + \beta\gamma\mu^2\theta]\tau^2 + \nonumber \\ + [\mu^2 (1 - \alpha) + \beta\gamma(\beta\gamma - \alpha\mu + \delta) + \alpha\delta\mu + (\beta\gamma - \alpha\mu + \delta)\theta] \tau \stackrel{\text{(\ref{cond_num_a1})-(\ref{cond_num_a2})}}{>} 0, 
\end{gather}

since the condition (\ref{stablecond}) implies
\begin{gather}
\label{cond_num_a1} \gamma > \frac{\alpha\mu - \delta}{\beta} \Rightarrow \beta\gamma - \alpha\mu + \delta > 0, \\
\label{cond_num_a2} \text{and} \ \gamma > \sqrt{\frac{( 2 - \alpha)\epsilon\mu + \alpha\epsilon \delta}{\beta\mu}} \Rightarrow \beta^2 \gamma^2 \mu^2 - 2 \beta \epsilon \mu^2 - \alpha \beta \delta \epsilon \mu + \alpha \beta \epsilon \mu^2 > 0. 
\end{gather}

Hence, the numerator of the ratio (\ref{fractiondetAAA1}) is strictly greater than 0 under the condition (\ref{stablecond}). \\

(b) Under the condition (\ref{stablecond2}), the numerator of the ratio (\ref{fractiondetAAA1}) is
\begin{gather}
(\beta\epsilon\tau - 1)^2 \mu^2 \tau + (\beta\epsilon\tau - 1)^2  + \nonumber \\
+ [(\beta^2\gamma^2\mu^2 - \alpha\beta\delta\epsilon\mu + \alpha\beta\epsilon\mu^2) + \delta\mu^2 \theta(1 - \alpha) + \beta\delta\gamma\mu^2(1 - \alpha)  + \beta\gamma\mu^2\theta]\tau^2 + \nonumber \\ + [( \beta^2\gamma^2 - \beta\gamma (\alpha\mu - \delta) -\alpha\mu^2) + \alpha\delta\mu + (\beta\gamma - \alpha\mu + \delta)\theta] \tau \stackrel{\text{(\ref{cond_num_b1})-(\ref{cond_num_b2})}}{>} 0, 
\end{gather}

since the condition (\ref{stablecond2}) implies
\begin{gather}
\gamma > \frac{\alpha \mu - \delta + \sqrt{(\alpha \mu - \delta)^2 + 4\alpha\mu^2}}{2 \beta} > \frac{\alpha\mu - \delta}{\beta} \nonumber \\
\label{cond_num_b1} \Rightarrow \beta^2\gamma^2 - \beta\gamma (\alpha\mu - \delta) -\alpha\mu^2 > 0 \ \text{and} \ \beta\gamma - \alpha\mu + \delta > 0 \\ 
\label{cond_num_b2} \text{and} \ \gamma > \sqrt{\max\left\{\frac{\alpha\epsilon(\delta - \mu)}{\beta\mu},0\right\}} \Rightarrow \beta^2\gamma^2\mu^2 - \alpha\beta\delta\epsilon\mu + \alpha\beta\epsilon\mu^2 > 0. 
\end{gather}

Hence, the numerator of the ratio (\ref{fractiondetAAA1}) is strictly greater than 0 under the condition (\ref{stablecond2}). \\

Therefore, $A_1 A_2$ has no negative eigenvalues under the condition either (\ref{stablecond}) or (\ref{stablecond2}). $\Box$ \\

\textit{Conclusion of the proof of Theorem \ref{thrm3}}. By Lemma \ref{lemmaa1}, $A_1$ is Hurwitz for all positive $\beta$, $\gamma$, $\mu$, $\epsilon$; $\delta \geq 0$, $\theta \geq 0$; and $\alpha \in [0,1)$. By Lemma \ref{a2cond}, $A_2$ is Hurwitz under the condition either (\ref{stablecond}) or (\ref{stablecond2}). It is easy to verify that the difference $A_1 - A_2$ has rank one. By Lemma \ref{lemmaa1a2noeig}, $A_1 A_2$ has no negative real eigenvalues under the condition either (\ref{stablecond}) or (\ref{stablecond2}). Hence, by Proposition \ref{prop4}, two LTI systems $u^\prime(t) = A_1 u(t)$ and $u^\prime(t) = A_2 u(t)$ have a CQLF. Therefore, the system (\ref{noboundary_system2}) is exponentially stable under the condition either (\ref{stablecond}) or (\ref{stablecond2}). $\Box$

\section{Numerical and simulation experiments and conjectures}\label{numerical}

In this section, we present some numerical and simulation experiments. These results are for both stylized and actual schemes, and all results are for the true system which includes boundary $X \geq 0$. We also put forward some conjectures based on these experiments. 

In all simulations, we always assume $r=1000$, but specify only the actual arrival rate $\Lambda=\lambda r$. On the plots labeled 'fluid', 
$X(t), Y(y), V(t)$ are replaced by their {\em fluid approximations} 
$$
X(t) = r x(t) + \frac{\lambda r(1-\alpha)}{\beta}, ~~ Y(t) = r y(t), ~~ V(t) = r v(t) + \frac{\lambda r}{\mu},
$$
respectively, where $(x(\cdot), y(\cdot), v(\cdot)$ is the corresponding fluid limit.

\subsection{Stylized scheme}\label{fluid_stylized}

\begin{figure}[h]
    \centering
    \begin{subfigure}[b]{0.45\textwidth}
        \includegraphics[width=\textwidth]{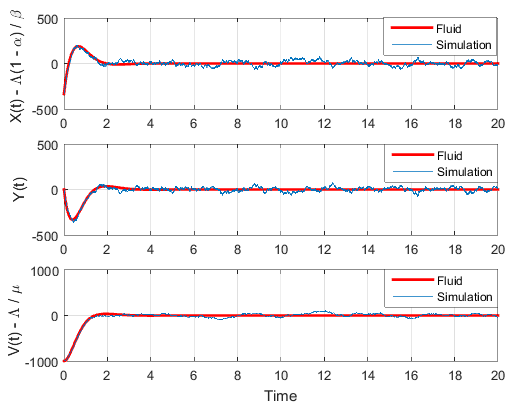}
        \caption{$(X(0),Y(0),Z(0)) = (0,0,0)$}
        \label{exp01a}
    \end{subfigure}
    \begin{subfigure}[b]{0.45\textwidth}
        \includegraphics[width=\textwidth]{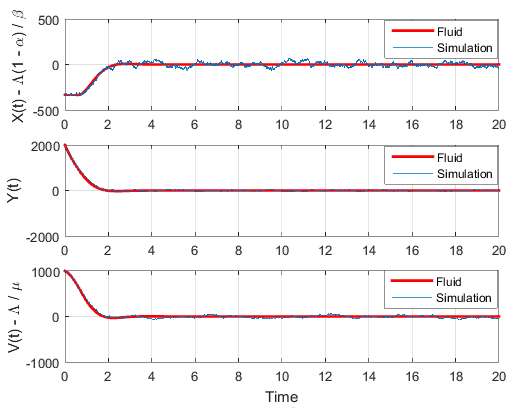}
        \caption{$(X(0),Y(0),Z(0)) = (0,2000,0)$}
        \label{exp01b}
    \end{subfigure}
    \begin{subfigure}[b]{0.45\textwidth}
        \includegraphics[width=\textwidth]{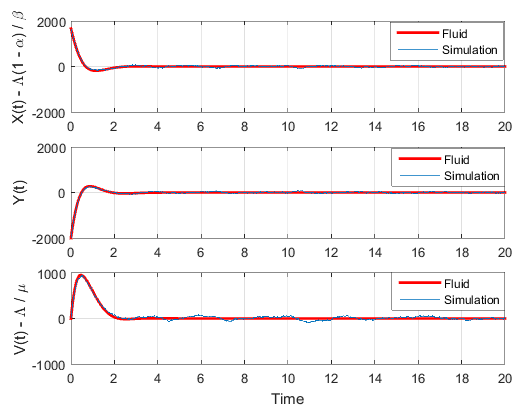}
        \caption{$(X(0),Y(0),Z(0)) = (2000,-2000,1000)$}
        \label{exp01c}
    \end{subfigure}
    \begin{subfigure}[b]{0.45\textwidth}
        \includegraphics[width=\textwidth]{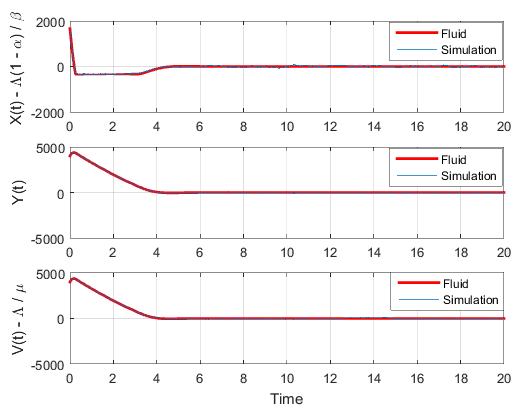}
        \caption{$(X(0),Y(0),Z(0)) = (2000,4000,1000)$}
        \label{exp01d}
    \end{subfigure}        
    \caption{Stylized scheme: Comparison of fluid approximations with simulations in Example \ref{exmp1}}\label{exp01}
\end{figure}

\begin{exmp}\label{exmp1}
Consider the following set of parameters, which satisfies condition (\ref{stablecond}):
\begin{gather*}
\Lambda = 2000 \ , \ \alpha = 0.5 \ , \ \beta = 3 \ , \ \mu = 2 \ , \ \gamma = 1 \ , \ \epsilon = 1.5 \ , \ \delta = 1 \ , \ \theta = 0.1 
\end{gather*}
\end{exmp}

with four initial conditions: (a) $(X(0),Y(0),Z(0)) = (0,0,0)$; (b) $(X(0),Y(0),Z(0)) = (0,2000,0)$; (c) $(X(0),Y(0),Z(0)) = (2000,-2000,1000)$; (d) $(X(0),Y(0),Z(0)) = (2000,4000,1000)$. The red line of the figure is the fluid approximation and the blue one is the simulation experiment. We see the converging trajectories on the Figure \ref{exp01}. Note that Figures \ref{exp01b} and \ref{exp01d} show that the trajectory hits the boundary on $X$. We also did the numerical/simulation experiments with many different sets of parameters satisfying the condition (\ref{stablecond}). All results, including those not shown on Figure \ref{exp01}, suggest the global stability of the system.

\begin{figure}[h]
    \centering
    \begin{subfigure}[b]{0.45\textwidth}
        \includegraphics[width=\textwidth]{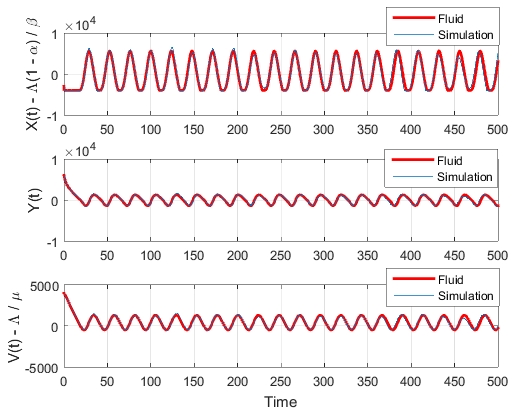}
        \caption{$\gamma = 1$}
        \label{exp02a}
    \end{subfigure}
    \begin{subfigure}[b]{0.45\textwidth}
        \includegraphics[width=\textwidth]{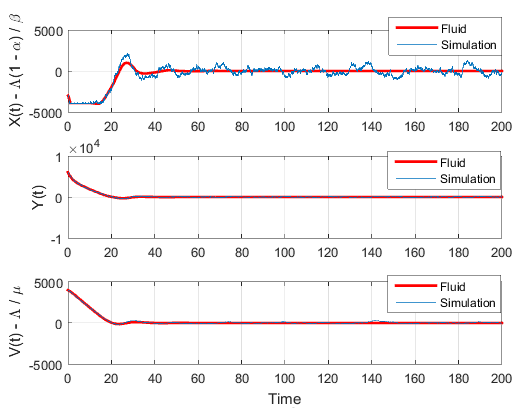}
        \caption{$\gamma = 5$}
        \label{exp02b}
    \end{subfigure}
    \begin{subfigure}[b]{0.45\textwidth}
        \includegraphics[width=\textwidth]{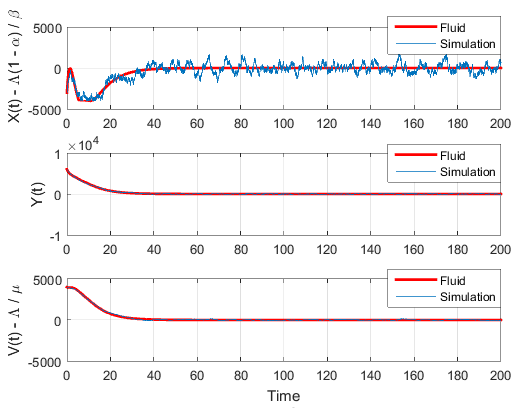}
        \caption{$\gamma = 10$}
        \label{exp02c}
    \end{subfigure}
    \begin{subfigure}[b]{0.45\textwidth}
        \includegraphics[width=\textwidth]{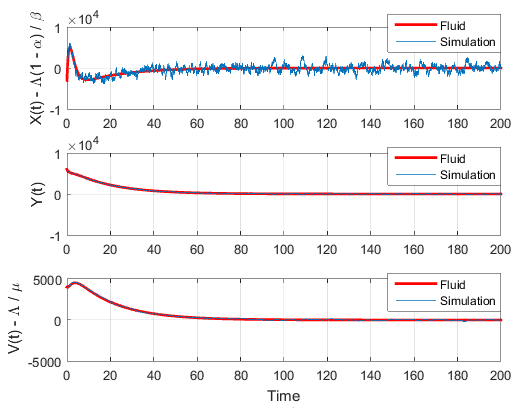}
        \caption{$\gamma = 20$}
        \label{exp02d}
    \end{subfigure}        
    \caption{Stylized scheme: Comparison of fluid approximations with simulations in Example \ref{exmp2}}\label{exp02}
\end{figure}

\begin{exmp}\label{exmp2}
We use sets of parameters:
\begin{gather*}
\Lambda = 2000 \ , \ \alpha = 0.9 \ , \ \beta = 0.05 \ , \ \mu = 0.5 \ , \ \epsilon = 1 \ , \ \delta = 0.01 \ , \ \theta = 0.01
\end{gather*}
\end{exmp}

with four different values of $\gamma$ ($\gamma_1 = 1$, $\gamma_2 = 5$, $\gamma_3 = 10$, and $\gamma_4 = 20$) (Figure \ref{exp02}). The sets of parameters with $\gamma_1 = 1$ and $\gamma_2 = 5$ do not satisfy the condition (\ref{stablecond}) while the sets of parameters with $\gamma_3 = 10$ and $\gamma_4 = 20$ satisfy the condition (\ref{stablecond}). We consider an initial condition $(X(0),Y(0),Z(0)) = (1000,6000,2000)$. On the Figures \ref{exp02b}, \ref{exp02c} and \ref{exp02d}, we see that the trajectories converge. However, Figure \ref{exp02a} shows the trajectory that never converges under the set of parameters with $\gamma_1 = 1$. \\

With many numerical/simulation experiments, the results, including those not shown on Figure \ref{exp02}, suggest both local and global stability of the system for all sufficiently large $\gamma$. \\

Our simulation experiments show that the fluid trajectory provides a very good approximation for the behavior of stylized scheme. 

\subsection{Actual scheme}\label{fluid_actual}

\begin{figure}[h]
    \centering
    \begin{subfigure}[b]{0.45\textwidth}
        \includegraphics[width=\textwidth]{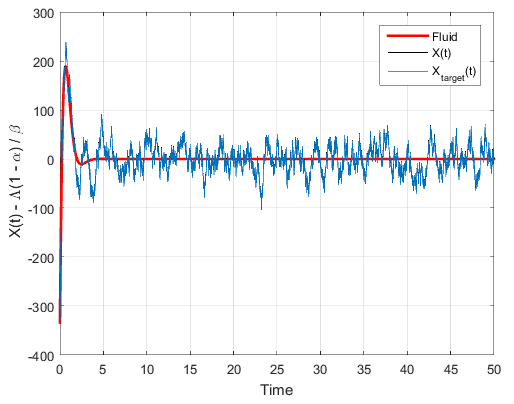}
        \caption{Fluid vs. $X(t)$ and $X_{target}(t)$}
        \label{exp03a1}
    \end{subfigure}
    \begin{subfigure}[b]{0.45\textwidth}
        \includegraphics[width=\textwidth]{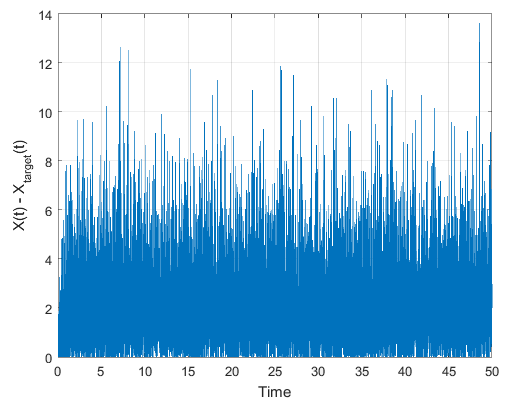}
        \caption{$X(t) - X_{target}(t)$}
        \label{exp03a2}
    \end{subfigure}        
    \caption{Actual scheme: $(X(0),Y(0),Z(0),X_{target}(0)) = (0,0,0,0)$}\label{exp03a}
\end{figure}

\begin{figure}[h]
    \centering
    \begin{subfigure}[b]{0.45\textwidth}
        \includegraphics[width=\textwidth]{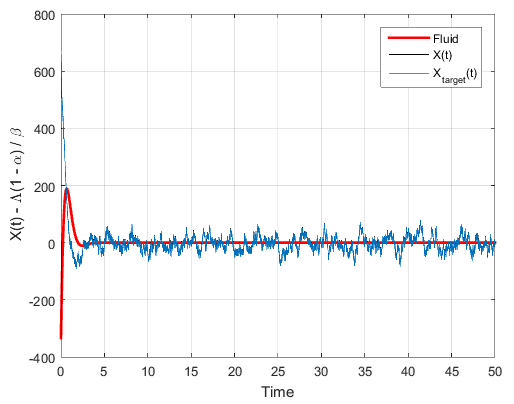}
        \caption{Fluid vs. $X(t)$ and $X_{target}(t)$}
        \label{exp03b1}
    \end{subfigure}
    \begin{subfigure}[b]{0.45\textwidth}
        \includegraphics[width=\textwidth]{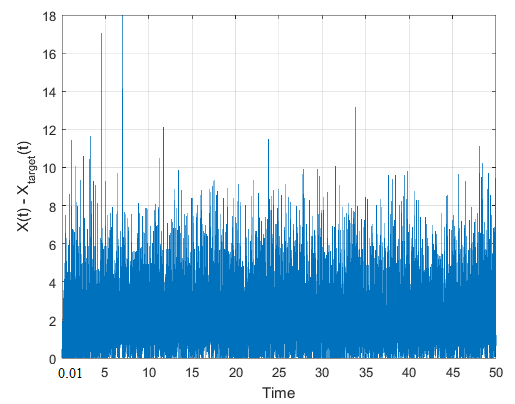}
        \caption{$X(t) - X_{target}(t)$}
        \label{exp03b2}
    \end{subfigure}        
    \caption{Actual scheme: $(X(0),Y(0),Z(0),X_{target}(0)) = (0,0,0,1000)$}\label{exp03b}
\end{figure}

\begin{exmp}\label{exmp3}
We conduct a simulation experiment for the actual scheme with the same set of parameters as in Example \ref{exmp1}: 
\begin{gather*}
\Lambda = 2000 \ , \ \alpha = 0.5 \ , \ \beta = 3 \ , \ \mu = 2 \ , \ \gamma = 1 \ , \ \epsilon = 1.5 \ , \ \delta = 1 \ , \ \theta = 0.1 
\end{gather*}
\end{exmp}

with two initial conditions $(X(0),Y(0),Z(0),X_{target}(0)) = (0,0,0,0)$ and $(X(0),Y(0),Z(0),X_{target}(0)) = (0,0,0,1000)$. (Note that this set of parameters satisfies the condition (\ref{stablecond}).) The results are shown in Figures \ref{exp03a} and \ref{exp03b}. We see that the magnitude of the difference between $X_{target}$ and the actual number of invited agents $X$ is very small (except at time 0) and can be negligible compared to their values. This explains why the trajectories of $X_{target}$ and $X$ are well approximated by the fluid trajectory, obtained for the stylized scheme.

For the stylized scheme, the results suggest the global stability of our system for all sufficiently large $\gamma$. However, the problem with large $\gamma$ is that the behavior of the stylized scheme may significantly deviate from the behavior of the actual scheme, as illustrated by the following example. 

\begin{figure}[h]
    \centering
    \begin{subfigure}[b]{0.45\textwidth}
        \includegraphics[width=\textwidth]{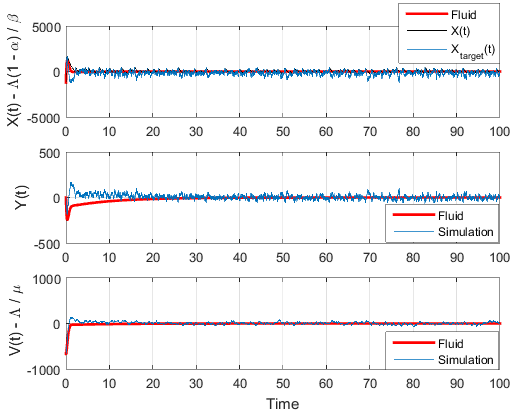}
        \caption{$\gamma = 10$}
        \label{exp04a}
    \end{subfigure}
    \begin{subfigure}[b]{0.45\textwidth}
        \includegraphics[width=\textwidth]{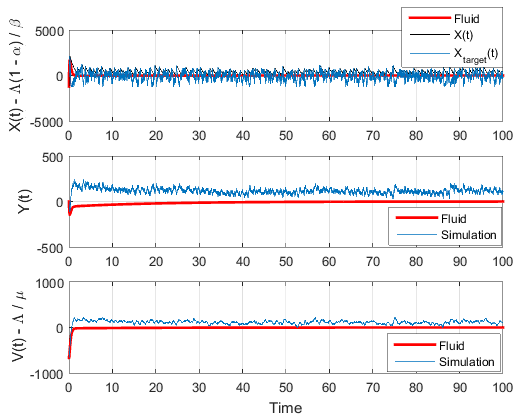}
        \caption{$\gamma = 20$}
        \label{exp04b}
    \end{subfigure}        
    \caption{Problem with large $\gamma$ of the actual scheme}\label{exp04}
\end{figure}

\begin{exmp}\label{exmp4}
Consider the following set of parameters:
\begin{gather*}
\Lambda = 2000 \ , \ \alpha = 0.7 \ , \ \beta = 0.5 \ , \ \mu = 3 \ , \ \epsilon = 1 \ , \ \delta = 1 \ , \ \theta = 2
\end{gather*}
\end{exmp}

with two values of $\gamma$ ($\gamma_1 = 10$, $\gamma_2 = 20$); and an initial condition $(X(0),Y(0),Z(0),X_{target}(0)) = (0,0,0,1000)$ (Figure \ref{exp04}). These results show that the behavior of the actual scheme deviates substantially from the behavior of the fluid trajectory with large $\gamma$. 

\begin{figure}[h]
    \centering
    \begin{subfigure}[b]{0.45\textwidth}
        \includegraphics[width=\textwidth]{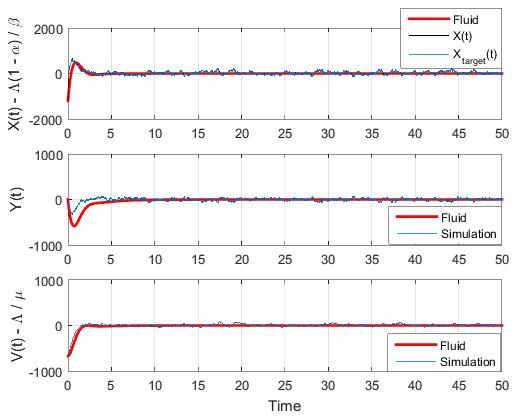}
        \caption{$\gamma = 2.3$}
        \label{exp04c1}
    \end{subfigure}
    \begin{subfigure}[b]{0.45\textwidth}
        \includegraphics[width=\textwidth]{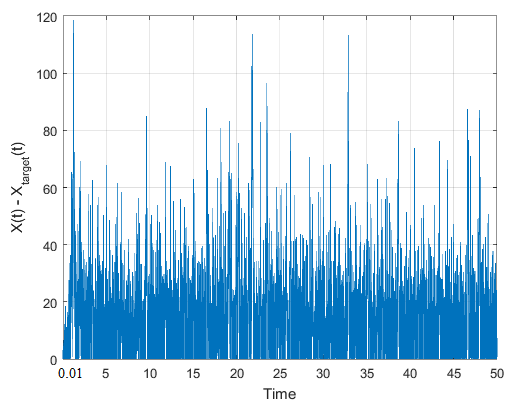}
        \caption{$X(t) - X_{target}(t)$}
        \label{exp04c2}
    \end{subfigure}        
    \caption{A ``good'' value of $\gamma$ for the actual scheme}\label{exp04c}
\end{figure}

Since $\alpha \mu > \delta$ and $\epsilon \leq \frac{(\alpha\mu - \delta)^2 \mu}{(2-\alpha)\mu\beta + \alpha\delta\beta}$, then we choose $\gamma = 2.3$ such that $\gamma > \frac{\alpha \mu - \delta}{\beta}$ (Corollary \ref{cor8}). We can see that, with a ``good'' value of $\gamma$, the behavior of the actual scheme deviates negligibly from the behavior of the fluid trajectory (Figure \ref{exp04c1}) and the difference between $X_{target}$ and $X$ is not large compared to their values (Figure \ref{exp04c2}). 

\subsection{Global vs. local stability of fluid limits}\label{global_local}

In this paper, we have derived some sufficient local stability conditions for the fluid limits. Based on a variety of simulation experiments above for the stylized scheme, we conjecture that local stability is sufficient for global stability of fluid limits for our model. In the next example, we compare the behavior of fluid limits for the system without boundary (given by (\ref{noboundary_system2})) with that of the system with boundary (given by (\ref{theorem2_system})).

\begin{figure}[h]
    \centering
    \begin{subfigure}[b]{0.45\textwidth}
        \includegraphics[width=\textwidth]{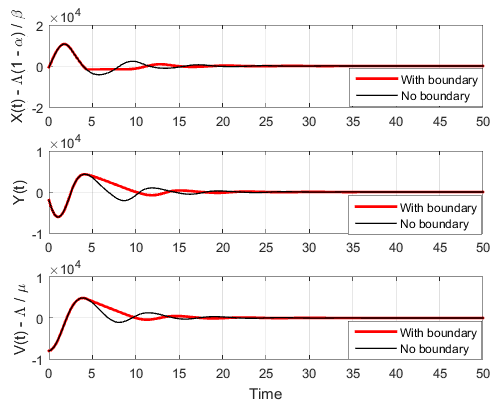}
        \caption{Set 1}
        \label{exp05a}
    \end{subfigure}
    \begin{subfigure}[b]{0.45\textwidth}
        \includegraphics[width=\textwidth]{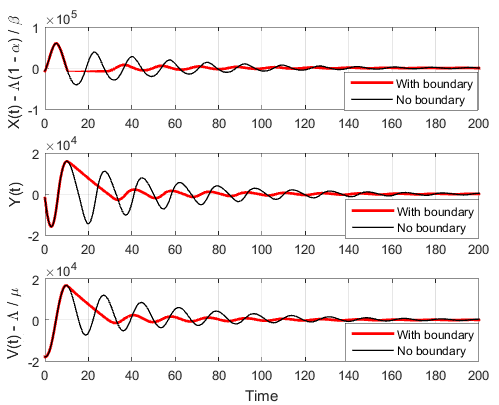}
        \caption{Set 2}
        \label{exp05b}
    \end{subfigure}        
    \caption{Fluid trajectories of the systems (\ref{theorem2_system}) and (\ref{noboundary_system2})}\label{exp05}
\end{figure}

\begin{exmp}\label{exmp5}
Consider two set of parameters, which satisfy the local stability conditions, so that the trajectory of the system (\ref{noboundary_system2}) converges to the equilibrium point $(0,0,0)$ (Figure \ref{exp05}). The red line of the figure is the trajectory of the system (\ref{theorem2_system}), which may hit the boundary $X=0$, and the black one is the trajectory of the system (\ref{noboundary_system2}), for which there is no boundary.  
\end{exmp}

With many experiments, the results, including those not shown in Figure \ref{exp05}, further suggest the global stability of the fluid limits, when the local stability holds. 

\subsection{Summary of conjectures, based on numerical and simulation experiments.}

\begin{conj}\label{conj01}
Our system is globally stable if it is locally stable. 
\end{conj}

\begin{conj}\label{conj02}
Given all other parameters are fixed, our system is globally stable for all sufficiently large $\gamma$. 
\end{conj}

Obviously, Conjecture \ref{conj01} is stronger than Conjecture \ref{conj02} because we have proved the local stability when $\gamma$ is large in this paper. We note again, however, that in a practical application the value of $\gamma$ should not be made too large, because the stylized scheme behavior, which we studied in this paper, may substantially deviate from the behavior of the actual scheme, where uninviting pending agents are not allowed. 

\section{Discussion and further work}\label{conclusion}

In this paper, we study a feedback-based agent invitation scheme for a model with randomly behaving agents and possible abandonment of customers and agents. This model is motivated by a variety of existing and emerging applications.
The focus of the paper is on the stability properties of the system fluid limits, arising as asymptotic limits of the system process, when the system scale (customer arrival rate) grows to infinity. The dynamic system, describing the behavior of fluid limit trajectories has a very complex structure -- it is a switched linear system, which in addition has a reflecting boundary. We derived some sufficient local stability conditions, using the machinery of switched linear systems and common quadratic Lyapunov functions. Our simulation and numerical experiments show good overall performance of the feedback scheme, when the local stability conditions hold. They also suggest that, for our model, the local stability is in fact sufficient for the global stability of fluid limits. Verifying these conjectures, as well as expanding the sufficient local stability conditions, is an interesting subject for future research.
 Further generalizations of the agent invitation model are also of interest from both theoretical and practical points of view.\\
 
 {\em Acknowledgement.} The authors would like to thank the referees for useful suggestions which helped to improve the exposition in the paper.

\bibliographystyle{abbrv}
\bibliography{reference}

\end{document}